\documentclass[12pt, twoside]{article}
\usepackage[pagewise]{lineno}
\usepackage{amsmath,amssymb}
\usepackage{float}
\usepackage{url}
\usepackage{multirow}
\usepackage{latexsym}
\usepackage{bm}
\usepackage{graphicx}
\graphicspath{{fig_eps/}}
\usepackage{epsfig}
\usepackage{cite}
\usepackage{epstopdf}
\usepackage{caption}
\usepackage{subfigure}
\usepackage[justification=centering]{caption}
\usepackage{textcomp,booktabs}
\usepackage[usenames,dvipsnames]{color}
\usepackage{colortbl}

\definecolor{mygray}{gray}{0.9}
\definecolor{mypink}{rgb}{0.99,0.91,0.95}
\definecolor{mycyan}{cmyk}{0.3,0,0,0}
\usepackage{booktabs}
\usepackage[top=2.5cm,bottom=2.5cm,left=2cm,right=2cm]{geometry}
\usepackage[colorlinks,linkcolor=red,anchorcolor=green,citecolor=blue,CJKbookmarks=True]{hyperref}
\newtheorem{theorem}{Theorem}[section]
\newtheorem{lemma}{Lemma}[section]

\newtheorem{definition}{Definition}[section]

\allowdisplaybreaks
\def\Proof{\noindent{\bf Proof.}~}
\def\qed{\hfill$\square$\smallskip}

\begin{document}
\title{\textbf{Traveling Wave in a Ratio-dependent Holling-Tanner System with Nonlocal Diffusion and Strong Allee Effect}}
\author{Hongliang Li\thanks{
Corresponding author.}, Min Zhao, Rong Yuan\\
Laboratory of Mathematics and Complex Systems (Ministry of Education)\\
School of Mathematical Sciences, Beijing Normal University \\
Beijing 100875, People's Republic of China\\
E-mail: lihl@mail.bnu.edu.cn}
\date{}
\maketitle
\begin{abstract}
In this paper, a ratio-dependent Holling-Tanner system with nonlocal diffusion is taken into account, where the prey is subject to a strong Allee effect. To be special, by applying Schauder's fixed point theorem and iterative technique, we provide a general theory on the existence of traveling waves for such system. Then appropriate upper and lower solutions and a novel sequence, similar to squeeze method, are constructed to demonstrate the existence of traveling waves for $c>c^*$. Moreover, the existence of traveling wave for $c=c^*$ is also established by spreading speed theory and comparison principle. Finally, the nonexistence of traveling waves for $c<c^*$ is investigated, and the minimal wave speed then is determined.

\textbf{Key words:} Holling-Tanner system; Traveling waves; Upper and lower solutions; Strong Allee effect.
\end{abstract}

\section{Introdcution}
\noindent

Predation is one of the most fundamental interactions in the natural world, which has become a ubiquitous fact. The classical Lotka-Volterra predator-prey system may not accurately describe the complex ecological since populations are limited by factors such as competition for resources, disease in reality. Additionally, predators are able to adapt their hunting behavior in response to changes in prey density and availability, and other ecological factors such as habitat fragmentation, climate change, and human activities can also influence predator-prey dynamics. Thus, to develop a more comprehensive understanding of complex ecological behavior and interactions, some extended and improved predator-prey systems have been studied extensively, seeing \cite{m89,m19,t97}.

Holling-Tanner system with Holling type II functional response has recently attracted increasing interest\cite{hh95, ag99}. As said in \cite{c97}, the ratio-dependent Holling-Tanner system provides a way to avoid the ``biological control paradox'' wherein classical prey-dependent exploitation models generally do not allow for a pest prey equilibrium density that is both low and stable. Moreover, based on strongly supported view in \cite{ag89, h97} that the functional responses over ecological timescales ought to depend on the densities of both prey and predator and more specifically on the ratio of prey to predator, Arditi and Ginzburg in \cite{ag89} proposed a ratio-dependent functional response
\begin{equation*}
  f\left(\displaystyle\frac{u}{v}\right)=\frac{m\left(\displaystyle\frac{u}{v}\right)}{\left(\displaystyle\frac{u}{v}\right)+a}=\frac{mu}{u+av},
\end{equation*}
and the corresponding Holling-Tanner system is studied in \cite{lp07}. Owing to spatially heterogeneous natural habitats, the reaction-diffusion system is established, traveling waves with special translation invariant form of various functional response Holling-Tanner system
are further studied to better characterize the population propagation. We refer the readers to \cite{cgy17, adp17, wf21, d12, zw22, zs18, zw20} for local diffusion system, and \cite{cy17, dlz19} for nonlocal diffusion system.

It should be noted that the prey is generally assumed to grow at a logistic pattern in most of the aforementioned Holling-Tanner systems, whereas, Allee effect has been studied extensively in recent years, because of its complexity and practicability, seeing \cite{ds07,csw14,am04,dps19}. Allee effect refers to reduced fitness or decline in population growth at low population densities. It is difficult to find mates and birth rates decrease  when population densities are low, which can lead to extinction\cite{al07}. It can be distinguished into two types: weak Allee effect and strong Allee effect \cite{wk01,wsw11}. The former refers to a scenario where the population has a positive and increasing growth when the size of the population is below a certain threshold value, and the latter pertains to a situation where growth is negative when the size of the population is below a certain threshold value.
Particularly, Zhao and Wu in \cite{zw20} investigated the existence of traveling waves of Holling-Tanner system with Lotka-Volterra functional response and strong Allee effect
\begin{equation}\label{eq:zw}
\left\{
\begin{array}{l}
u_t=d_{1}\Delta u+u\left(1-u\right)\left(\displaystyle\frac{u}{b}-1\right)-auv, \\[0.2cm]
v_t=d_{2}\Delta v+s v\left(1-\displaystyle\frac{v}{u}\right).
\end{array}
\right.
\end{equation}
As far as the author knows, there are almost no related results on the existence of traveling waves for ratio-dependent Holling-Tanner systems with nonlocal diffusion and strong Allee effect. Hence, in this work, we consider the following system
\begin{equation}\label{eq:rd}
\left\{
\begin{array}{l}
u_t=d_{1}\mathcal{N}_1[u](x,t)+u\left(1-u\right)\left(\displaystyle\frac{u}{b}-1\right)-\displaystyle \frac{muv}{u+av}, \\[0.4cm]
v_t=d_{2}\mathcal{N}_2[v](x,t)+s v\left(1-\displaystyle\frac{v}{u}\right).
\end{array}
\right.
\end{equation}
where $u(x,t)$ and $v(x,t)$ stand for the population densities of the prey and the predator, respectively. All the parameters in system \eqref{eq:rd} are positive. $d_1$ and $d_2$ are the diffusion coefficients of $u(x,t)$ and $v(x,t)$; $b\in(0,1)$ represents Allee threshold value; $m$ is a measure of the quality of the prey as food for the predator; $a$ and $s$ denote the saturation rate of the predator and a measure of the growth rate of the predator. Moreover, $\mathcal{N}_i[w](x,t)$, $i=1,2,$ formulates the spatial nonlocal diffusion of individuals
\begin{align*}
  \mathcal{N}_i[w](x,t)=\int_{\mathbb R} J_i(x-y)w(y,t)dy-w(x,t),
\end{align*}
in which the kernel function $J_i: \mathbb{R}\rightarrow\mathbb{R}$, $i = 1, 2$, satisfies the following assumptions:
\begin{description}
  \item \emph{(J1)}\ $J_i\in C^1(\mathbb R)$, $J_i(x)=J_i(-x)\geq0$ and $\displaystyle\int_{\mathbb{R}}J_i(x)dx=1$.
  \item \emph{(J2)}\ $J_i$ satisfy the decay bounds:
  \begin{equation*}
    \int_{\mathbb{R}}J_i(x)e^{\lambda x}dx<\infty \text{ for any } \lambda\in(0,\lambda_0) \text{ and } \lim\limits_{\lambda\rightarrow\lambda_0}\int_{\mathbb{R}}J_i(x)e^{\lambda x}dx=\infty
  \end{equation*}
   for some $\lambda_0\in(0,+\infty]$ and $\displaystyle\int_{\mathbb{R}}|J'_i(x)|dx<\infty$.
\end{description}

Moreover, if one chooses the kernel function as $J_i(x)=\delta(x)+\delta''(x)$, where $\delta(x)$ is the Dirac delta function\cite{m03}, then system \eqref{eq:rd} can be reduced to the following system \cite{zl22}:
\begin{equation}\label{eq:st}
\left\{
\begin{array}{l}
u_t=d_{1}\Delta u+u\left(1-u\right)\left(\displaystyle\frac{u}{b}-1\right)-\displaystyle \frac{muv}{u+av}, \\[0.4cm]
v_t=d_{2}\Delta v+s v\left(1-\displaystyle\frac{v}{u}\right).
\end{array}
\right.
\end{equation}
As for the corresponding kinetic system, that is, system \eqref{eq:rd} or \eqref{eq:st} without diffusion
\begin{equation}\label{eq:ck}
\left\{
\begin{array}{l}
u'(t)=u(1-u)\left(\displaystyle\frac{u}{b}-1\right)-\displaystyle\frac{muv}{u+av}, \\[0.4cm]
v'(t)=s v\left(1-\displaystyle\frac{v}{u}\right).
\end{array}
\right.
\end{equation}
Obviously, $(b, 0)$ and $(1, 0)$ are nonnegative equilibria of system \eqref{eq:ck}. Denote
\begin{equation*}
  b_1=1+\frac{2 m}{1+a}-2 \sqrt{\frac{m}{1+a}\left(1+\frac{m}{1+a}\right)}.
\end{equation*}
It is easy to verify that $b_1\in(0,1)$ and system \eqref{eq:ck} has a unique positive equilibrium $(u^*,u^*)$ if $b=b_1$, two positive equilibria $(u_{1}^{*},u_{1}^{*})$ and $(u_{2}^{*},u_{2}^{*})$ with $u_{1}^{*}<u_{2}^{*}<1$ if $0<b<b_1$ and no positive equilibrium if $b_1<b<1$, where
\begin{align*}
  &u^*=\frac{1+b_1}{2},\\
  &u_{1}^{*}=\frac{1}{2}\left(b+1-\sqrt{b^{2}-2\left(1+\frac{2 m}{1+a}\right) b+1}\right),\\
  &u_{2}^{*}=\frac{1}{2}\left(b+1+\sqrt{b^{2}-2\left(1+\frac{2 m}{1+a}\right) b+1}\right).
\end{align*}

In this paper, we mainly establish the existence of traveling waves connecting predator-free state and coexistence state of system \eqref{eq:rd}. In this paper, we assume $0<b<b_1$, which is equivalent to $4mb<(1-b)^2(1+a)$ owing to $b\in(0,1)$. According to \cite{zl22}, at this moment, positive equilibrium $(u_1^*, u_1^*)$ is unstable. Thus, we should focus on the existence of traveling waves connecting $(1,0)$ and $(u_2^*, u_2^*)$. A positive solution is called a traveling wave if it has the form
\begin{equation*}
  (u,v)(x,t)=(\phi, \psi)(\xi), \ \xi=x+ct,
\end{equation*}
where $c>0$ is a wave speed. Then $(\phi,\psi)(\xi)$ should satisfy
\begin{equation}\label{eq:tr}
\left\{\begin{split}
&c\phi'(\xi)=d_{1} \mathcal{N}_1[\phi](\xi)+\phi(\xi)f(\phi,\psi)(\xi),\\[0.4cm]
&c\psi'(\xi)=d_{2} \mathcal{N}_2[\psi](\xi)+\psi(\xi)g(\phi,\psi)(\xi),
\end{split}
\right.
\end{equation}
where
\begin{align*}
  \mathcal{N}_i\left[w\right](\xi)=\int_{\mathbb R} J_i(\xi-y)w(y)dy-w(\xi), \ i=1,2,
\end{align*}
and
\begin{equation*}
  f(\phi,\psi)=\left(1-\phi\right)\left(\displaystyle\frac{\phi}{b}-1\right)-\displaystyle\frac{m\psi}{\phi+a\psi} \ \text{   and  }\ g(\phi,\psi)=s \left(1-\displaystyle\frac{\psi}{\phi}\right).
\end{equation*}
If $(\phi,\psi)(\xi)$ further meets the boundary conditions
\begin{equation}\label{eq:b}
  (\phi,\psi)(-\infty)=(1,0)\ \text{ and }\ (\phi,\psi)(+\infty)=(u_2^*,u_2^*),
\end{equation}
then it is called a traveling wave connecting $(1,0)$ and $(u_2^*, u_2^*)$, also named a invasion wave \cite{d83}. As a result, demonstrating the existence of traveling waves of system \eqref{eq:rd} is exactly the same as demonstrating the existence of positive solution of system \eqref{eq:tr}-\eqref{eq:b}. Our main results are stated in the following.
\begin{theorem}\label{th:ma}
Assume that (J1)-(J2) hold. If $b\in(0,1)$ and
\begin{equation}\label{eq:as1}
  m<\min\left\{\frac{(1-b)^2(1+b)}{8b},\ \frac{(1-b)^2(1+a)}{4b}, \ \frac{(1+a)^3}{8}\left(\sqrt{b^2+4\left(\frac{1-b}{1+a}\right)^2}-b\right)\right\},
\end{equation}
then system \eqref{eq:rd} has traveling waves connecting $(1,0)$ and $(u_2^*,u_2^*)$ for $c\geq c^*$.
Additionally, system \eqref{eq:rd} has no traveling waves connecting $(1,0)$ and $(u_2^*,u_2^*)$ for $0<c< c^*$.
\end{theorem}

It should be emphasized that our method can be also applied to Holling-Tanner system with local diffusion and strong Allee effect,
and the related theory on spreading speed and comparison principle can be find in \cite{a78, ql11}. In particular, for system \eqref{eq:zw}, the existence of traveling wave with $c=c^*$ is proved through a limiting argument similar to Theorem 4.1 in \cite{h11}, and this process is omitted owing to its complexity. Here we use an entirely different approach with \cite{zw20} to deal with it, and then present a complete proof.

This paper is organized as follows. In Section 2, we provide a general theory on the existence of positive solution for system \eqref{eq:tr}. In Section 3, we firstly establish the existence of a positive solution of system \eqref{eq:tr}-\eqref{eq:b} by constructing appropriate upper and lower solutions and a novel sequence for $c>c^*$, where $c^*$ is denoted in \eqref{eq:cx}. In the last section, we show the case $c=c^*$ also holds based on spreading speed theory and comparison principle. Moreover, we also obtain the nonexistence of positive solution of system \eqref{eq:tr}-\eqref{eq:b}, which determines the minimal wave speed.

\section{A general result}
\noindent

In this section, by using Schauder's fixed point theorem and the upper and lower solutions method, we  transform the problem of finding a positive solution of system \eqref{eq:tr} into existence problem of a pair of upper and lower solutions. Firstly, we define
\begin{equation*}
  \begin{array}{l}
  X_{b}=\left\{(\phi,\psi)\in C\left(\mathbb{R}, \mathbb{R}^{2}\right) : (1+b)/2\leq \phi\leq 1 \text{ and } 0\leq \psi\leq 1\right\}.
  \end{array}
\end{equation*}
and introduce the definition of upper and lower solutions.
\begin{definition}\label{Def}
Function pairs $\left(\overline{\phi},\overline{\psi}\right)(\xi)$ and $\left(\underline{\phi},\underline{\psi}\right)(\xi)$ in $X_b$ are upper solution and lower solution of  system \eqref{eq:tr} if they satisfy
\begin{description}
  \item(a)\  $\underline{\phi}(\xi)\leq\overline{\phi}(\xi)$, $\underline{\psi}(\xi)\leq\overline{\psi}(\xi)$ for $\xi\in\mathbb R$.
  \item(b)\ There is a finite set $E=\left\{\xi_i: 1\leq i\leq m\right\}$ such that for $\xi\in\mathbb R\backslash E$
  \begin{align}
    &d_1
    \mathcal{N}_1\left[\overline{\phi}\right](\xi)-c\overline{\phi}'(\xi)+\overline{\phi}(\xi)f\left(\overline{\phi},\underline{\psi}\right)(\xi)\leq0,\label{eq:s1}\\[0.2cm]
    &d_1
    \mathcal{N}_1\left[\underline{\phi}\right](\xi)-c\underline{\phi}'(\xi)+\underline{\phi}(\xi)f\left(\underline{\phi},\overline{\psi}\right)(\xi)\geq0,\label{eq:s2}\\[0.2cm] &d_2 \mathcal{N}_2\left[\overline{\psi}\right](\xi)-c\overline{\psi}'(\xi)+\overline{\psi}(\xi)g\left(\overline{\phi},\overline{\psi}\right)(\xi)\leq0,\label{eq:s3}\\[0.2cm]
    &d_2
    \mathcal{N}_2\left[\underline{\psi}\right](\xi)-c\underline{\psi}'(\xi)+\underline{\psi}(\xi)g\left(\underline{\phi},\underline{\psi}\right)(\xi)\geq0.\label{eq:s4}
  \end{align}
\end{description}
\end{definition}

Next, we consider the nonlinear operators $F_1$ and $F_2$ defined on $X_b$ by
\begin{align*}
  F_ {1}(\phi,\psi)(\xi):&=\beta\phi(\xi)+d_ {1}\mathcal{N}_ {1}[\phi](\xi)+\phi(\xi)f\left(\phi,\psi\right)(\xi),\\[0.2cm]
  F_ {2}(\phi,\psi)(\xi):&=\beta\psi(\xi)+d_ {2}\mathcal{N}_ {2}[\psi](\xi)+ \psi(\xi)g\left(\phi,\psi\right)(\xi),
\end{align*}
where positive constant $\beta$ is large enough, such that $F_1(\phi,\psi)$ is nondecreasing in $\phi$ and is decreasing in $\psi$, and $F_2(\phi,\psi)$ is nondecreasing in $\phi$ and $\psi$ for $(\phi,\psi)(\xi)\in X_b$, respectively.
We also define the following operators
\begin{align*}
  P_ {1}(\phi,\psi)(\xi)&=\frac{1}{c}\int_{-\infty}^{\xi}e^{\frac{\beta(y-\xi)}{c}}F_ {1}(\phi,\psi)(y) dy, \ \xi\in\mathbb R,\\[0.2cm]
  P_ {2}(\phi,\psi)(\xi)&=\frac{1}{c}\int_{-\infty}^{\xi}e^{\frac{\beta(y-\xi)}{c}}F_ {2}(\phi,\psi)(y) dy,\ \xi\in\mathbb R.
\end{align*}
Set $P = (P_1,P_2)$, then $P : X_b\rightarrow  C\left(\mathbb{R}, \mathbb{R}^{2}\right)$ and it is easy to check that $(\hat{\phi},\hat{\psi})=P(\hat{\phi},\hat{\psi})$ solves
\begin{equation*}
\left\{\begin{split}
&c\hat{\phi}'(\xi)=-\beta\hat{\phi}(\xi)+F_ {1}(\hat{\phi},\hat{\psi})(\xi), \\[0.2cm]
&c\hat{\psi}'(\xi)=-\beta\hat{\psi}(\xi)+F_ {2}(\hat{\phi},\hat{\psi})(\xi).
\end{split}
\right.
\end{equation*}
Hence a fixed point of $P$ is a solution of system \eqref{eq:tr}. Therefore, it remains to show that $P$ has a fixed point in $X_b$. To do it, we choose constant $\mu\in(0,\beta/c)$ and define space
\begin{equation*}
  B_\mu(\mathbb{R}, \mathbb{R}^2):=\left\{(\phi,\psi) \in X_b : |(\phi,\psi)|_\mu=\sup\limits_{\xi\in\mathbb{R}}\left\{ \max\left(|\phi(\xi)|, |\psi(\xi)|\right) \right\}e^ {-\mu|\xi|}< \infty\right\}.
\end{equation*}
Then $\left(B_\mu(\mathbb{R}, \mathbb{R}^2), |\cdot|_\mu\right)$ is a Banach space from \cite{m01}. By applying Schauder's fixed point theorem, we will seek a positive solution of system \eqref{eq:tr} in set
\begin{equation*}
  \Sigma=\left\{(\phi,\psi)\in X_b: \underline{\phi}\leq \phi\leq\overline{\phi},\ \underline{\psi}\leq \psi\leq\overline{\psi}\right\},
\end{equation*}
which is a non-empty convex, closed and bounded set in $(B_\mu(\mathbb{R}, \mathbb{R}^2), |\cdot|_\mu)$. \\

Firstly, we show that $P(\Sigma)\subseteq\Sigma$. Let $(\phi,\psi)\in\Sigma$, since $(\phi,\psi)$, $\left(\overline{\phi},\overline{\psi}\right)$ and $\left(\underline{\phi},\underline{\psi}\right)$ all belong to $X_b$, we have
\begin{equation*}
  F_1(\overline{\phi},\underline{\psi})(\xi) \geq F_1(\phi,\psi)(\xi)\geq F_1(\underline{\phi},\overline{\psi})(\xi)
\end{equation*}
 by the choice of $\beta$. We assert that
\begin{equation*}
  \underline{\phi}(\xi)\leq P_1(\phi,\psi)(\xi)\leq\overline{\phi}(\xi).
\end{equation*}
Assume that $-\infty:=\xi_{m+1}<\xi_m<\xi_{m-1}<\cdots<\xi_1<\xi_0:=+\infty$ in set $E$ of Definition \ref{Def}. Through integration by parts formula, we drive for $0\leq i\leq m$
\begin{equation*}
  \frac{1}{c}\int_{\xi_{i+1}}^{\xi_i}e^{\frac{\beta}{c}y}\left(c\phi'(y)+\beta\phi(y) \right) dy=e^{\frac{\beta}{c}\xi_i}\phi(\xi_i)-e^{\frac{\beta}{c}\xi_{i+1}}\phi(\xi_{i+1}).
\end{equation*}
If $\xi\in(\xi_{k+1},\xi_k)$ for $0\leq k\leq m$, from the definitions of upper and lower solutions, we arrive at
\begin{equation*}
\begin{split}
  P_1(\phi,\psi)(\xi)\geq P_ {1}(\underline{\phi},\overline{\psi})(\xi)
  &=\frac{1}{c}\int_{-\infty}^{\xi}e^{\frac{\beta(y-\xi)}{c}}F_1(\underline{\phi},\overline{\psi})(y) dy \\[0.2cm]
  &\geq\frac{1}{c}\int_{-\infty}^{\xi}e^{\frac{\beta(y-\xi)}{c}}\left(c\underline{\phi}'(y)+\beta\underline{\phi}(y) \right) dy\\[0.2cm]
  &=\frac{1}{c}\left(\sum_{j=k+1}^{m}\int_{\xi_{j+1}}^{\xi_j}+\int_{\xi_{k+1}}^{\xi}\right)e^{\frac{\beta(y-\xi)}{c}}\left(c\underline{\phi}'(y)+\beta\underline{\phi}(y) \right) dy\\[0.2cm]
  &=\underline{\phi}(\xi)-e^{\frac{\beta(\xi_{m+1}-\xi)}{c}}\underline{\phi}(\xi_{m+1})=\underline{\phi}(\xi)
\end{split}
\end{equation*}
owing to $\xi_{m+1}=-\infty$. Similarly, we get
\begin{equation*}
\begin{split}
  P_1(\phi,\psi)(\xi)\leq P_ {1}(\overline{\phi},\underline{\psi})(\xi)
  &=\frac{1}{c}\int_{-\infty}^{\xi}e^{\frac{\beta(y-\xi)}{c}}F_1(\overline{\phi},\underline{\psi})(y) dy \\[0.2cm]
  &\leq\frac{1}{c}\int_{-\infty}^{\xi}e^{\frac{\beta(y-\xi)}{c}}\left(c\overline{\phi}'(y)+\beta\overline{\phi}(y) \right) dy\\[0.2cm]
  &=\frac{1}{c}\left(\sum_{j=k+1}^{m}\int_{\xi_{j+1}}^{\xi_j}+\int_{\xi_{k+1}}^{\xi}\right)e^{\frac{\beta(y-\xi)}{c}}\left(c\overline{\phi}'(y)+\beta\overline{\phi}(y) \right) dy\\[0.2cm]
  &=\overline{\phi}(\xi)-e^{\frac{\beta(\xi_{m+1}-\xi)}{c}}\overline{\phi}(\xi_{m+1})=\overline{\phi}(\xi)
\end{split}
\end{equation*}
owing to $\xi_{m+1}=-\infty$. Then in a similar way, we can verify that
\begin{equation*}
  \underline{\psi}(\xi)\leq P_2(\phi,\psi)(\xi)\leq\overline{\psi}(\xi).
\end{equation*}
Therefore, we have $P(\Sigma)\subseteq\Sigma$.
Next, by the choice of $\mu$, we show that $P$ is completely continuous with respect to the norm $|\cdot|_\mu$. First of all, we show the continuity of $P$ on $\Sigma$ with respect to the norm $|\cdot|_\mu$.
Let $\Phi_1=(\phi_1,\psi_1)$ and $\Phi_2=(\phi_2,\psi_2)$ be in $\Sigma$, then a direct calculation yields that
\begin{align*}
  &|F_ {1}(\phi_1,\psi_1)(\xi)-F_ {1}(\phi_2,\psi_2)(\xi)|\\[0.2cm]
  =&|\beta\phi_1(\xi)+d_ {1}\mathcal{N}_ {1}[\phi_1](\xi)+\phi_1(\xi)f(\phi_1,\psi_1)(\xi)-\beta\phi_2(\xi)-d_ {1}\mathcal{N}_ {1}[\phi_2](\xi)-\phi_2(\xi)f(\phi_2,\psi_2)(\xi)|\\[0.2cm]
  \leq&\left(\beta+d_1+C\right)|\phi_1(\xi)-\phi_2(\xi)|+C|\psi_1(\xi)-\psi_2(\xi)|+d_1\int_{\mathbb R}J_1(\xi-s)|\phi_1(s)-\phi_2(s)|ds,
\end{align*}
where
\begin{equation*}
  C=\max\limits_{(\phi,\psi)\in X_b}\left\{f(\phi,\psi)+\phi\frac{\partial f(\phi,\psi)}{\partial \phi},\ \phi\frac{\partial f(\phi,\psi)}{\partial \psi}\right\}.
\end{equation*}
Furthermore, we have
\begin{align*}
  &|P_1(\phi_1,\psi_1)(\xi)-P_1(\phi_2,\psi_2)(\xi)|e^ {-\mu|\xi|}\\[0.2cm]
  =&\frac{1}{c}\Big |\int_{-\infty}^{\xi}e^{\frac{\beta(y-\xi)}{c}}\left(F_ {1}(\phi_1,\psi_1)(y)-F_ {1}(\phi_2,\psi_2)(y)\right) dy\Big |e^ {-\mu|\xi|}\\[0.2cm]
  \leq&\frac{e^ {-\mu|\xi|}}{c}\int_{-\infty}^{\xi}e^{\frac{\beta(y-\xi)}{c}} |F_ {1}(\phi_1,\psi_1)(y)-F_ {1}(\phi_2,\psi_2)(y) | dy\\[0.2cm]
  \leq&\frac{\left(\beta+d_1+C\right)e^ {-\mu|\xi|}}{c}\int_{-\infty}^{\xi}e^{\frac{\beta(y-\xi)}{c}} |\phi_1(y)-\phi_2(y)|dy\\[0.2cm]
  &+\frac{Ce^ {-\mu|\xi|}}{c}\int_{-\infty}^{\xi}e^{\frac{\beta(y-\xi)}{c}} |\psi_1(y)-\psi_2(y)|dy\\[0.2cm]
  &+\frac{d_1e^ {-\mu|\xi|}}{c}\int_{-\infty}^{\xi}e^{\frac{\beta(y-\xi)}{c}} \left(\int_{\mathbb R}J_1(y-s)|\phi_1(s)-\phi_2(s)|ds\right) dy
  =I_1+I_2+I_3.
\end{align*}
Now, let us continue estimating $I_1$, $I_2$ and $I_3$, respectively. Note that $|y|-|\xi|\leq\xi-y$ for all $y\leq\xi$, we have
\begin{equation}\label{eq:z1}
  \int_{-\infty}^{\xi}e^{\frac{\beta(y-\xi)}{c}} e^{\mu|y|}e^ {-\mu|\xi|}dy\leq\int_{-\infty}^{\xi}e^{\frac{(\beta-c\mu)(y-\xi)}{c}}dy=\frac{c}{\beta-c\mu}.
\end{equation}
Recalling that $\Phi_1=(\phi_1,\psi_1)$ and $\Phi_2=(\phi_2,\psi_2)$, for the part $I_1$, we utilize \eqref{eq:z1} to get
\begin{align*}
  I_1&=\frac{\left(\beta+d_1+C\right)e^ {-\mu|\xi|}}{c}\int_{-\infty}^{\xi}e^{\frac{\beta(y-\xi)}{c}} |\phi_1(y)-\phi_2(y)|dy\\[0.2cm]
  &=\frac{\left(\beta+d_1+C\right)e^ {-\mu|\xi|}}{c}\int_{-\infty}^{\xi}e^{\frac{\beta(y-\xi)}{c}} e^ {\mu|y|}|\phi_1(y)-\phi_2(y)|e^ {-\mu|y|}dy\\[0.2cm]
  &\leq\frac{|\Phi_1-\Phi_2|_{\mu}\left(\beta+d_1+C\right)}{c}\int_{-\infty}^{\xi}e^{\frac{\beta(y-\xi)}{c}} e^{\mu|y|}e^ {-\mu|\xi|}dy\\[0.2cm]
  &\leq\frac{\beta+d_1+C}{\beta-c\mu}|\Phi_1-\Phi_2|_{\mu}.
\end{align*}
For the part $I_2$, similarly, we get
\begin{align*}
  I_2\leq\frac{C}{\beta-c\mu}|\Phi_1-\Phi_2|_{\mu}.
\end{align*}
For the part $I_3$, we also have
\begin{align*}
  I_3&=\frac{d_1e^ {-\mu|\xi|}}{c}\int_{-\infty}^{\xi}e^{\frac{\beta(y-\xi)}{c}} \left(\int_{\mathbb R}J_1(y-s)|\phi_1(s)-\phi_2(s)|ds\right) dy\\[0.2cm]
  &\leq\frac{d_1e^ {-\mu|\xi|}|\Phi_1-\Phi_2|_{\mu}}{c}\int_{-\infty}^{\xi}e^{\frac{\beta(y-\xi)}{c}} \left(\int_{\mathbb R}J_1(z)e^{\mu|y-z|}dz\right) dy\\[0.2cm]
  &\leq\frac{d_1|\Phi_1-\Phi_2|_{\mu}}{c}\left(\int_{\mathbb R}J_1(z)e^{\mu|z|}dz\right)\int_{-\infty}^{\xi}e^{\frac{\beta(y-\xi)}{c}}e^{\mu |y|}e^ {-\mu|\xi|} dy\\[0.2cm]
  &\leq \frac{d_1}{\beta-c\mu}\left(\int_{\mathbb R}J_1(z)e^{\mu|z|}dz\right)|\Phi_1-\Phi_2|_{\mu}<\infty
\end{align*}
owing to \emph{(J2)}. The above three estimates suggest that
\begin{equation*}
  |P_1(\phi_1,\psi_1)(\xi)-P_1(\phi_2,\psi_2)(\xi)|e^ {-\mu|\xi|}\leq M_1|\Phi_1-\Phi_2|_{\mu},
\end{equation*}
where
\begin{equation*}
  M_1=\frac{\beta+d_1+2C+d_1\displaystyle\int_{\mathbb R}J_1(z)e^{\mu|z|}dz}{\beta-c\mu},
\end{equation*}
Similarly, there is a constant $M_2$ such that
\begin{equation*}
  |P_2(\phi_1,\psi_1)(\xi)-P_2(\phi_2,\psi_2)(\xi)|e^ {-\mu|\xi|}\leq M_2|\Phi_1-\Phi_2|_{\mu}.
\end{equation*}
Hence $P$ is continuous on $\Sigma$ with respect to the norm $|\cdot|_\mu$. Simultaneously, for $(\phi,\psi)\in\Sigma$,
\begin{align*}
  &\Big|\frac{d }{d \xi} P_1(\phi,\psi)(\xi)\Big|=\Big|-\frac{\beta}{c}P_1(\phi,\psi)(\xi)+\frac{1}{c}F_1(\phi,\psi)(\xi)\Big|   \leq\frac{\beta}{c}\overline{\phi}(\xi)+\frac{1}{c}F_1(1,0)\leq\frac{2\beta}{c},\\[0.2cm]
  &\Big|\frac{d }{d \xi}P_2(\phi,\psi)(\xi)\Big|=\Big|-\frac{\beta}{c}P_2(\phi,\psi)(\xi)+\frac{1}{c}F_2(\phi,\psi)(\xi)\Big|\leq\frac{\beta}{c}\overline{\psi}(\xi)+\frac{1}{c}F_2(1,1)\leq\frac{2\beta}{c}.
\end{align*}
Thus, $P(\Sigma)$ is equicontinuous. From Proposition 2.12 in \cite{ez86}, an argument similar to Lemma 3.4 in \cite{h03} can verify that $P:\Sigma\rightarrow\Sigma$ is compact. Hence, $P$ is completely continuous with respect to the norm $|\cdot|_\mu$. Now, applying Schauder's fixed point theorem, there exists $(\hat{\phi},\hat{\psi})\in\Sigma$ such that $P(\hat{\phi},\hat{\psi})=(\hat{\phi},\hat{\psi})$,
which implies that $(\hat{\phi},\hat{\psi})$ is a fixed point of $P$ in $\Sigma$. Therefore, we conclude the following lemma.
\begin{lemma}\label{le:ssm}
Assume that system \eqref{eq:tr} has a pair of upper and lower solutions $\left(\overline{\phi},\overline{\psi}\right)(\xi)$ and $\left(\underline{\phi},\underline{\psi}\right)(\xi)$ satisfying
\begin{equation*}
    \underline{\phi}(\xi_i^-)\leq\underline{\phi}(\xi_i^+),\ \overline{\phi}(\xi_i^+)\leq\overline{\phi}(\xi_i^-), \ \underline{\psi}(\xi_i^-)\leq\underline{\psi}(\xi_i^+),\ \overline{\psi}(\xi_i^+)\leq\overline{\psi}(\xi_i^-)\ \text{for}\ \xi_i\in E.
\end{equation*}
Then it has a solution $(\phi,\psi)(\xi)$ satisfying for $\xi\in\mathbb R$
\begin{equation*}
  \underline{\phi}(\xi)\leq \phi(\xi)\leq\overline{\phi}(\xi),\quad \underline{\psi}(\xi)\leq \psi(\xi)\leq\overline{\psi}(\xi).
\end{equation*}
\end{lemma}
\section{The existence of traveling waves for $c>c^*$}
\noindent

The major goal of this section is to demonstrate the existence of traveling waves with $c>c^*$ of system \eqref{eq:rd}, where $c^*$ is denoted in \eqref{eq:cx}. Firstly, a pair of appropriate upper and lower solutions is constructed. Then a positive solution $(\phi,\psi)(\xi)$ of system \eqref{eq:tr} is immediately obtained. Inspired by squeeze method \cite{d13,adp17}, we construct a novel sequence to study asymptotic behavior of such positive solution at $\xi=+\infty$, which further leads to the existence of traveling waves for $c>c^*$.

\subsection{Upper and lower solutions}
\noindent

To find a pair of upper and lower solutions, we define
\begin{equation*}
  \Delta(\lambda,c)=d_2\left(\int_{\mathbb R} J_2(y) e^{-\lambda y}dy-1\right)-c\lambda+s,
\end{equation*}
By direct calculation, it is easy to get
\begin{align*}
\begin{array}{ll}
&\Delta(0, c)=s>0 \text{ and } \lim\limits_{\lambda \rightarrow+\infty} \Delta(\lambda, c)=+\infty \text { for all } c, \\[0.3cm]
&\displaystyle\frac{\partial \Delta(\lambda, c)}{\partial c}=-\lambda<0 \text { and } \lim\limits_{c \rightarrow+\infty} \Delta(\lambda, c)=-\infty \text { for } \lambda>0, \\[0.4cm]
&\left.\displaystyle\frac{\partial \Delta(\lambda, c)}{\partial \lambda}\right|_{\lambda=0}=-c<0 \text{ for all } c>0, \\[0.4cm]
&\displaystyle\frac{\partial^{2} \Delta(\lambda, c)}{\partial \lambda^{2}}=d_2 \int_{\mathbb{R}} J_2(y) y^{2} e^{-\lambda y} dy>0 \text { for all } \lambda \text { and } c .
\end{array}
\end{align*}
Thus, we have the following properties for $\Delta(\lambda,c)$.
\begin{lemma}\label{le:cx}
There exists a positive constant
\begin{equation}\label{eq:cx}
  c^*=\inf\limits_{\lambda>0}\left\{\frac{1}{\lambda}\left[d_2\left(\int_{\mathbb R} J_2(y) e^{-\lambda y}dy-1\right)+s\right]\right\}
\end{equation}
such that the following assertions hold.
\begin{description}
  \item (a) If  $0<c<c^{*}$, then  $\Delta(\lambda, c)>0$ for $\lambda>0$.
  \item (b) If  $c>c^{*}$, then $\Delta(\lambda, c)=0$  has two positive real roots  $\lambda_1<\lambda_2$, and $\Delta(\cdot, c)<0$ in $(\lambda_1,\lambda_2)$ and $\Delta(\cdot, c)>0$ in $(0,\lambda_1)\cup(\lambda_2,+\infty)$.
  \item (c) If  $c=c^{*}$, then there exists  $\lambda^{*}>0$  such that  $\Delta\left(\lambda^{*}, c^{*}\right)=0$ and  $\Delta\left(\lambda, c^{*}\right)>0$  for  $\lambda \neq \lambda^{*}$.
\end{description}
\end{lemma}

Additionally, we define
\begin{equation*}
  \Pi(\lambda, c)=d_{1} \left(\int_{\mathbb{R}} J_{1}(y)\mathrm{e}^{-\lambda y} \mathrm{d}y-1\right)-c \lambda,
\end{equation*}
then one can easily check that for $c>c^*$
\begin{equation*}
  \Pi(0, c)=0 \text{ and } \left.\displaystyle\frac{\partial \Pi(\lambda, c)}{\partial \lambda}\right|_{\lambda=0}=-c<0.
\end{equation*}
Thus, we choose $\eta\in(0,\lambda_1)$ to satisfy $\Pi(\eta, c)<0$.\\

Based on the above discussions, we now introduce function pairs $(\overline{\phi},\overline{\psi})$ and $(\underline{\phi},\underline{\psi})$,
\begin{align}\label{eq:sup-sub1}
\begin{split}
  \overline{\phi}(\xi)&=1,
  \qquad\qquad \qquad\quad \
  \underline{\phi}(\xi)=\left\{
  \begin{array}{ll}
  \displaystyle\frac{1+b}{2}, &\xi>0,\\[0.2cm]
  1-\displaystyle\frac{1-b}{2}e^{\eta\xi}, &\xi\leq 0,
  \end{array}
  \right.
  \\[0.4cm]
  \overline{\psi}(\xi)&=\left\{
  \begin{array}{ll}
  1, &\xi> 0,\\[0.2cm]
  e^{\lambda_1\xi}, &\xi\leq 0,
  \end{array}
  \right.
  \quad
  \underline{\psi}(\xi)=\left\{
  \begin{array}{ll}
  0, &\xi>\xi_1:=-(1/\varepsilon)\ln r,\\[0.2cm]
  e^{\lambda_1\xi}\left(1-re^{\varepsilon \xi}\right), &\xi\leq \xi_1,
  \end{array}
  \right.
\end{split}
\end{align}
where $\varepsilon$ and $r$ satisfy
\begin{equation*}
  0<\varepsilon<\min\left\{\lambda_1,\ \lambda_2-\lambda_1\right\} \text{ and }
  r>\max\left\{1,\ \frac{-2s}{(1+b)\Delta(\lambda_1+\varepsilon,c)}\right\}.
\end{equation*}

\begin{lemma}\label{le:s1}
For $c>c^*$, $(\overline{\phi},\overline{\psi})$ and $(\underline{\phi},\underline{\psi})$ satisfy \eqref{eq:s1} and \eqref{eq:s2} under the assumption \eqref{eq:as1}.
\end{lemma}
\Proof Recalling that
\begin{equation*}
  f(\phi,\psi)=\left(1-\phi\right)\left(\displaystyle\frac{\phi}{b}-1\right)-\displaystyle\frac{m\psi}{\phi+a\psi},
\end{equation*}
since $\overline{\phi}(\xi)=1$ for  $\xi\in\mathbb{R}$, clearly, \eqref{eq:s1} holds since
\begin{equation*}
  d_1
    \mathcal{N}_1\left[\overline{\phi}\right](\xi)-c\overline{\phi}'(\xi)+\overline{\phi}(\xi)f\left(\overline{\phi},\underline{\psi}\right)(\xi)=-\displaystyle\frac{m\underline{\psi}(\xi)}{1+a\underline{\psi}(\xi)}\leq0.
\end{equation*}
For \eqref{eq:s2}, since $\underline{\phi}(\xi)$ is non-increasing in $\mathbb R$, then
\begin{align*}
  \int_{\mathbb R} J_1(\xi-y)\underline{\phi}(y)dy\geq&\frac{1+b}{2}\int_{\mathbb R} J_1(\xi-y)dy=\frac{1+b}{2},
\end{align*}
and
\begin{align*}
  \int_{\mathbb R} J_1(\xi-y)\underline{\phi}(y)dy
  \geq\int_{\mathbb R} J_1(\xi-y)\left(1-\displaystyle\frac{1-b}{2}e^{\eta y}\right)dy
  =1-\frac{1-b}{2}e^{\eta \xi}\int_{\mathbb R} J_1(y)e^{-\eta y}dy.
\end{align*}
Hence we have
\begin{align}\label{eq:n1}
  \int_{\mathbb R} J_1(\xi-y)\underline{\phi}(y)dy
  \geq\max\left\{\frac{1+b}{2},\ 1-\frac{1-b}{2}e^{\eta \xi}\int_{\mathbb R} J_1(y)e^{-\eta y}dy\right\}.
\end{align}
If $\xi>0$, then $\underline{\phi}(\xi)=(1+b)/2$ and $\overline{\psi}(\xi)=1$. It follows from \eqref{eq:n1} that
\begin{equation*}
   \mathcal{N}_1\left[\underline{\phi}\right](\xi)=\int_{\mathbb R} J_i(\xi-y)\underline{\phi}(y)dy-\underline{\phi}(\xi)\geq0,
\end{equation*}
by the assumption \eqref{eq:as1}, we further get
\begin{align*}
  d_1
    \mathcal{N}_1\left[\underline{\phi}\right](\xi)-c\underline{\phi}'(\xi)+\underline{\phi}(\xi)f\left(\underline{\phi},\overline{\psi}\right)(\xi)
  \geq\underline{\phi}\left(\frac{(1-b)^2}{4b}-\frac{2m}{1+b+2a}\right)>0.
\end{align*}
If $\xi<0$, then $\underline{\phi}=1-(1-b)e^{\eta \xi}/2$ and $\overline{\psi}=e^{\lambda_1\xi}$.
It follows from \eqref{eq:n1} that
\begin{align*}
  \mathcal{N}_1\left[\underline{\phi}\right](\xi)=&\int_{\mathbb R} J_i(\xi-y)\underline{\phi}(y)dy-\underline{\phi}(\xi)\\[0.2cm]
  \geq&1-\frac{1-b}{2}e^{\eta \xi}\int_{\mathbb R} J_1(y)e^{-\eta y}dy-\left(1-\displaystyle\frac{1-b}{2}e^{\eta \xi}\right)\\[0.2cm]
  =&\frac{1-b}{2}e^{\eta \xi}\left(1-\int_{\mathbb R} J_1(y)e^{-\eta y}dy\right).
\end{align*}
Thus, we have
\begin{align*}
  &d_1
    \mathcal{N}_1\left[\underline{\phi}\right](\xi)-c\underline{\phi}'(\xi)+\underline{\phi}(\xi)f\left(\underline{\phi},\overline{\psi}\right)(\xi)\\[0.2cm]
  \geq&\frac{1-b}{2}e^{\eta\xi}\left[d_1\left(1-\int_{\mathbb R}J_1(y)e^{-\eta y}dy\right)+c\eta\right]\\[0.2cm]
  &+\underline{\phi}e^{\eta\xi}\left[\frac{1-b}{2}\left(\frac{1}{b}-\frac{1-b}{2b}e^{\eta\xi}-1\right)-\frac{me^{(\lambda_1-\eta)\xi}}{\underline{\phi}(\xi)+ae^{\lambda_1\xi}}\right]\\[0.2cm]
  =&-\frac{1-b}{2}e^{\eta\xi}\Pi(\eta,c)+\underline{\phi}e^{\eta\xi}K_1
\end{align*}
where
\begin{align*}
  K_1&=\frac{1-b}{2}\left(\frac{1}{b}-\frac{1-b}{2b}e^{\eta\xi}-1\right)-\frac{me^{(\lambda_1-\eta)\xi}}{\underline{\phi}(\xi)+ae^{\lambda_1\xi}}.
\end{align*}
Since $0<\eta<\lambda_1$, then $(\lambda_1-\eta)\xi<0$. Owing to $\underline{\phi}\geq(1+b)/2$ over $\mathbb R$, we get
\begin{align*}
  K_1\geq\frac{(1-b)^2}{4b}-\frac{2m}{1+b}>0
\end{align*}
by the assumption \eqref{eq:as1}. On the other hand, $\Pi(\eta,c)<0$ by the choice of $\eta$, which ensures that \eqref{eq:s2} holds. Therefore, we complete the proof.

\qed

\begin{lemma}\label{le:s2}
For $c>c^*$, $(\overline{\phi},\overline{\psi})$ and $(\underline{\phi},\underline{\psi})$ satisfy \eqref{eq:s3} and \eqref{eq:s4}.
\end{lemma}
\Proof For \eqref{eq:s3}, it is easy to check that
\begin{equation*}
  \int_{\mathbb R} J_1(\xi-y)\overline{\psi}(y)dy \leq \min \left\{1,e^{\lambda_{1} \xi} \int_{\mathbb{R}} J_{2}(y) e^{-\lambda_{1} y} d y\right\}.
\end{equation*}
If $\xi>0$, then $\overline{\psi}(\xi)=\overline{\phi}(\xi)=1$. Recall that $g(\phi,\psi)=s\left(1-\psi/\phi\right)$,
clearly,
\begin{align*}
  &d_2 \mathcal{N}_2\left[\overline{\psi}\right](\xi)-c\overline{\psi}'(\xi)+\overline{\psi}(\xi)g\left(\overline{\phi},\overline{\psi}\right)(\xi)\leq0.
\end{align*}
If $\xi<0$, then $\overline{\psi}(\xi)=e^{\lambda_1\xi}$ and $\overline{\phi}(\xi)=1$. Note that $\Delta(\lambda_1,c)=0$ for $c>c^*$, we arrive at
\begin{align*}
  &d_2 \mathcal{N}_2\left[\overline{\psi}\right](\xi)-c\overline{\psi}'(\xi)+\overline{\psi}(\xi)g\left(\overline{\phi},\overline{\psi}\right)(\xi)\\[0.2cm]
  \leq&e^{\lambda_1\xi}\left[d_2\left(\int_{\mathbb R}J_2(y) e^{-\lambda_1 y}dy-1\right)-c\lambda_1+s\right]-se^{2\lambda_1\xi}\\[0.2cm]
  =&e^{\lambda_1\xi}\Delta(\lambda_1,c)-se^{2\lambda_1\xi}=-se^{2\lambda_1\xi}\leq0.
\end{align*}
Hence \eqref{eq:s3} holds. For \eqref{eq:s4}, it is easy to check that
\begin{equation*}
  \int_{\mathbb R} J_1(\xi-y)\underline{\psi}(y)dy\geq\max\left\{0,\int_{\mathbb R}J_2(y)e^{\lambda_1(\xi-y)}\left(1-re^{\varepsilon (\xi-y)}\right)dy\right\}.
\end{equation*}
If $\xi>\xi_1$, $\underline{\psi}(\xi)=0$, clearly, \eqref{eq:s4} hold. If $\xi<\xi_1<0$,
\begin{align*}
  \underline{\psi}(\xi)=e^{\lambda_1\xi}\left(1-re^{\varepsilon \xi}\right) \text{  and  } \underline{\psi}'(\xi)=e^{\lambda_1\xi}\left(\lambda_1-r(\lambda_1+\varepsilon)e^{\varepsilon\xi}\right).
\end{align*}
Since $\varepsilon<\lambda_2-\lambda_1$, we have $\Delta(\lambda_1+\varepsilon,c)<0$ for $c>c^*$.
Owing to $\underline{\phi}\geq(1+b)/2$ over $\mathbb R$, we arrive at
\begin{align*}
  &d_2
    \mathcal{N}_2\left[\underline{\psi}\right](\xi)-c\underline{\psi}'(\xi)+\underline{\psi}(\xi)g\left(\underline{\phi},\underline{\psi}\right)(\xi)\\[0.2cm]
  \geq&d_2 \left(\int_{\mathbb R}J_2(y)e^{\lambda_1(\xi-y)}\left(1-re^{\varepsilon (\xi-y)}\right)dy-e^{\lambda_1\xi}\left(1-re^{\varepsilon \xi}\right)\right)\\[0.2cm]
  &-ce^{\lambda_1\xi}\left(\lambda_1-r(\lambda_1+\varepsilon)e^{\varepsilon\xi}\right)+se^{\lambda_1\xi}\left(1-re^{\varepsilon \xi}\right)-\frac{se^{2\lambda_1\xi}\left(1-re^{\varepsilon \xi}\right)^2}{\underline{\phi}(\xi)}\\[0.2cm]
  \geq&e^{\lambda_1\xi}\left[d_2\left(\int_{\mathbb R} J_2(y) e^{-\lambda_1 y}dy-1\right)-c\lambda_1+s\right]-\frac{2se^{2\lambda_1\xi}}{1+b}\\[0.2cm]
  &-re^{(\lambda_1+\varepsilon) \xi}\left[d_2\left(\int_{\mathbb R} J_2(y) e^{-(\lambda_1+\varepsilon) y}dy-1\right)-c(\lambda_1+\varepsilon)+s\right]\\[0.2cm]
  =&e^{\lambda_1\xi}\Delta(\lambda_1,c)-re^{(\lambda_1+\varepsilon) \xi}\Delta(\lambda_1+\varepsilon,c)-\frac{2se^{2\lambda_1\xi}}{1+b}\\[0.2cm]
  =&e^{(\lambda_1+\varepsilon)\xi}\left(-r\Delta(\lambda_1+\varepsilon,c)-\frac{2se^{(\lambda_1-\varepsilon)\xi}}{1+b}\right)\\[0.2cm]
  \geq&e^{(\lambda_1+\varepsilon)\xi}\left(-r\Delta(\lambda_1+\varepsilon,c)-\frac{2s}{1+b}\right)>0
\end{align*}
by the choice of $r$. Hence \eqref{eq:s4} holds. Therefore, we complete the proof.

\qed

\subsection{Asymptotic behavior}
\noindent

In the light of upper and lower solutions, a positive solution $(\phi,\psi)(\xi)$ of system \eqref{eq:tr} is immediately obtained with the aid of Lemmas \ref{le:ssm}. In this subsection, we mainly investigate the asymptotic behavior of the positive solution as $\xi=\infty$.
We firstly give the existence result on the positive solution.
\begin{theorem}\label{th:str}
Under the assumption \eqref{eq:as1}, system \eqref{eq:tr} has a positive solution $(\phi,\psi)(\xi)$ satisfying $(1+b)/2<\phi(\xi)<1$ and $0<\psi(\xi)<1$ over $\mathbb R$ for all $c>c^*$.
\end{theorem}
\Proof From Lemma \ref{le:ssm}, system \eqref{eq:tr} has a positive solution $(\phi,\psi)(\xi)$ satisfying $\underline{\phi}(\xi)\leq \phi(\xi)\leq\overline{\phi}(\xi)$  and  $\underline{\psi}(\xi)\leq \psi(\xi)\leq\overline{\psi}(\xi)$ over $\mathbb R$. We first show that $\phi(\xi)>(1+b)/2$ over $\mathbb R$. For contradiction, assume that there is a $\xi_0\in\mathbb R$ such that $\phi(\xi_0)=(1+b)/2$, we have $\phi'(\xi_0)=0$ owing to $\underline{\phi}(\xi)\geq(1+b)/2$ over $\mathbb R$. By $\phi$-equation, since $\overline{\psi}(\xi)\leq1$,  we get
\begin{equation*}
  0\geq-d_{1} \mathcal{N}_1[\phi](\xi_0)=\phi(\xi_0)f\left(\phi,\psi\right)(\xi_0)
  \geq\phi(\xi_0)f\left((1+b)/2,1\right)>0
\end{equation*}
under the assumption \eqref{eq:as1}. Thus, we have $\phi(\xi)>(1+b)/2$ over $\mathbb R$. Below, we show that $\psi(\xi_0)>0$ over $\mathbb R$. For contradiction, assume that there is a $\xi_0\in\mathbb R$ such that $\psi(\xi_0)=0$, then $\psi'(\xi_0)=0$ owing to $\underline{\psi}(\xi)\geq0$ over $\mathbb R$. By $\psi$-equation, we get
\begin{equation*}
  \int_{\mathbb R}J_2(\xi_0-y)\psi(y)dy=0 \text{ for all } y\in\mathbb R.
\end{equation*}
Thus, we have $\psi(\xi)\equiv0$ over $\mathbb R$, which contradicts $\psi(\xi)>\underline{\psi}(\xi)>0$ for $\xi<\xi_1$. Next, we show that $\phi(\xi)<1$ over $\mathbb R$. For contradiction, assume that there is a $\xi_0\in\mathbb R$ such that $\phi(\xi_0)=1$, then $\phi'(\xi_0)=0$ owing to $\phi(\xi)\leq1$ over $\mathbb R$. By $\phi$-equation, we get
\begin{align*}
  0\leq -d_{1} \mathcal{N}_1[\phi](\xi_0)&=-\displaystyle\frac{m\psi(\xi_0)}{1+a\psi(\xi_0)}<0.
\end{align*}
Thus, we have $\phi(\xi)<1$ over $\mathbb R$. Moreover, the case $\psi(\xi)<1$ can be treated similarly. Therefore, we complete the proof.

\qed

As for the positive solution $(\phi,\psi)(\xi)$ obtained in the above theorem, it follows from \eqref{eq:sup-sub1} that
\begin{align*}
&1=\lim_{\xi\rightarrow-\infty}\underline{\phi}(\xi)\leq\liminf_{\xi\rightarrow-\infty}\phi(\xi)\leq\limsup_{\xi\rightarrow-\infty}\phi(\xi)\leq\lim_{\xi\rightarrow-\infty}\overline{\phi}(\xi)=1,\\[0.2cm]
&0=\lim_{\xi\rightarrow-\infty}\underline{\psi}(\xi)\leq\liminf_{\xi\rightarrow-\infty}\psi(\xi)\leq\limsup_{\xi\rightarrow-\infty}\psi(\xi)\leq\lim_{\xi\rightarrow-\infty}\overline{\psi}(\xi)=0.
\end{align*}
Thus, we have $(\phi,\psi)(-\infty)=(1,0)$. Next, our main goal is to prove
\begin{equation*}
  \lim_{\xi\rightarrow+\infty}\phi(\xi)=\lim_{\xi\rightarrow+\infty}\psi(\xi)=u_2^*.
\end{equation*}
To do it,  we define
\begin{align*}
&\phi_+=\limsup \limits_{\xi \rightarrow+\infty} \phi(\xi),\quad \phi_-=\liminf\limits_{\xi \rightarrow+\infty} \phi(\xi), \\[0.2cm]
&\psi_{+}=\limsup \limits_{\xi \rightarrow+\infty} \psi(\xi),\quad \psi_-=\liminf\limits_{\xi \rightarrow+\infty} \psi(\xi) .
\end{align*}

\begin{lemma}\label{le:k1}
 $\phi_{-}\leq\psi_{-}\leq\psi_{+}\leq\phi_+$.
\end{lemma}
\textbf{Proof} We first show that $\phi_{-}\leq\psi_{-}$. For contradiction, we assume that $\phi_{-}>\psi_{-}$. First of all, we consider the case that $\psi(\xi)$ is eventually monotone, then
$\psi(+\infty)$ exists and
\begin{equation*}
  \int_{0}^{+\infty} \psi'(s)ds=\psi(+\infty)-\psi(0)<\infty
\end{equation*}
since $\psi(\xi)$ is bounded on $\mathbb R$. Furthermore, we have
$\limsup_{\xi\rightarrow+\infty}\psi'(\xi)=0$ if  $\psi'(\xi)\leq0$  for $\xi\gg 1$ or
$\liminf_{\xi\rightarrow+\infty}\psi'(\xi)=0$ if  $\psi'(\xi)\geq0$  for $\xi\gg 1$.
Hence we always can find a
sequence $\left\{\xi_n\right\}$ with $\xi_n\rightarrow+\infty$ as $n\rightarrow+\infty$ such that
\begin{equation}\label{eq:t1}
  \lim\limits_{n\rightarrow+\infty}\psi(\xi_n)=\psi_-=\psi_+<\phi_- \text{ and } \lim\limits_{n\rightarrow+\infty}\psi'(\xi_n)=0.
\end{equation}
Integrating $\psi$-equation of system \eqref{eq:tr} from 0 to $\xi_n$, we obtain
\begin{equation}\label{eq:k}
  c[\psi(\xi_n)-\psi(0)]-d_2\int_{0}^{\xi_n}
    \mathcal{N}_2\left[\psi\right](\xi)d\xi=\int_{0}^{\xi_n}\psi(\xi)g(\phi,\psi)(\xi)d\xi.
\end{equation}
Note that
\begin{equation}\label{eq:j1}
\begin{split}
   \int_{0}^{\xi_n} \mathcal{N}_2\left[\psi\right](\xi)d\xi&=\int_{0}^{\xi_n}\int_{\mathbb R}J_2(y)[\psi(\xi-y)-\psi(\xi)]dyd\xi\\[0.2cm]
  &=\int_{0}^{\xi_n}\int_{\mathbb R}J_2(y)(-y)\int_{0}^{1}\psi'(\xi-\tau y)d\tau dyd\xi\\[0.2cm]
  &=\int_{\mathbb R}J_2(y)(-y)\int_{0}^{1}\int_{0}^{\xi_n}\psi'(\xi-\tau y)d\xi d\tau dy\\[0.2cm]
  &=\int_{\mathbb R}J_2(y)(-y)\int_{0}^{1}[\psi(\xi_n-\tau y)-\psi(-\tau y)]d\tau dy.
\end{split}
\end{equation}
Since $0<\psi(\xi)<1$ over $\mathbb R$, from \emph{(J1)} and \emph{(J2)}, we further arrive at
\begin{equation}\label{eq:j2}
\begin{split}
&\int_{\mathbb{R}} J_2(y)(-y)\int_{0}^{1}[\psi(\xi_n-\tau y)-\psi(-\tau y)]d\tau dy\\[0.2cm]
\leq&2\int_{\mathbb{R}} J_2(y)|y|dy\leq4\int_{\mathbb{R}} J_2(y)e^{y}dy<\infty.
\end{split}
\end{equation}
One the other hand, by \eqref{eq:t1}, we have
\begin{equation*}
  \liminf\limits_{n\rightarrow +\infty}g(\phi,\psi)(\xi_n)\geq g(\phi_-,\psi_-)=s\left(1-\frac{\psi_- }{\phi_-}\right)>0.
\end{equation*}
Hence the left-hand of \eqref{eq:k} is bounded, whereas the right-hand of \eqref{eq:k} is unbounded, which leads to a contradiction. Next, we consider the case that $\psi(\xi)$ is oscillatory as $\xi\rightarrow +\infty$. We then can find a
sequence $\left\{\xi_n\right\}$ of the minimal points of $\psi(\xi)$ with $\xi_n\rightarrow+\infty$ as $n\rightarrow+\infty$ such that
\begin{equation}\label{eq:t2}
  \lim\limits_{n\rightarrow+\infty}\psi(\xi_n)=\psi_-<\phi_- \text{ and } \psi'(\xi_n)=0.
\end{equation}
And we also have
\begin{equation*}
  \liminf\limits_{n\rightarrow +\infty}g(\phi,\psi)(\xi_n)\geq g(\phi_-,\psi_-)=s\left(1-\frac{\psi_- }{\phi_-}\right)>0,
\end{equation*}
a repeating argument as above get what we desired. For the inequality $\psi_+\leq\phi_+$, we assume that $\psi_+>\phi_+$ by contradiction.
Similar to \eqref{eq:t1} and \eqref{eq:t2}, one always can choose sequence $\left\{\xi_n\right\}$  with $\xi_n\rightarrow+\infty$ as $n\rightarrow+\infty$ such that
\begin{equation*}
  \lim\limits_{n\rightarrow+\infty}\psi(\xi_n)=\psi_+>\phi_+ \text{ and } \lim\limits_{n\rightarrow+\infty}\psi'(\xi_n)=0.
\end{equation*}
At this moment, \eqref{eq:k} still holds, whereas
\begin{equation*}
  \limsup\limits_{n\rightarrow +\infty}g\left(\phi,\psi\right)(\xi_n)\leq g(\phi_+,\psi_+)=s\left(1-\frac{\psi_+}{\phi_+}\right)<0.
\end{equation*}
Hence a contradiction arises either $\psi(\xi)$ is eventually monotone or oscillatory as $\xi\rightarrow +\infty$.
Therefore, we complete the proof.

\qed

\begin{lemma}\label{le:0}
Under the assumption \eqref{eq:as1}, we have $(1+b)/2<\phi_{-}\leq\phi_{+}<1$.
\end{lemma}
\textbf{Proof}
From Theorem \ref{th:str}, we have
$$(1+b)/2\leq\phi_{-}\leq\phi_{+}\leq1.$$
Firstly, we show that $\phi_{-}>(1+b)/2$. For contradiction, we assume that $\phi_{-}=(1+b)/2$. If $\phi(\xi)$ is eventually monotone, we have $\phi(+\infty)$ exists and
\begin{equation*}
  \int_{0}^{+\infty} \phi'(s)ds=\phi(+\infty)-\phi(0)<\infty
\end{equation*}
since $\phi$ is bounded on $\mathbb{R}$.
Furthermore, we have
$\limsup_{\xi\rightarrow+\infty}\phi'(\xi)=0$ if  $\phi'(\xi)\leq0$  for $\xi\gg 1$ or
$\liminf_{\xi\rightarrow+\infty}\phi'(\xi)=0$ if  $\phi'(\xi)\geq0$  for $\xi\gg 1$.
Then we can find a sequence $\left\{\xi_n\right\}$ with ${\xi_n\rightarrow+\infty}$ as $n\rightarrow+\infty$ such that
\begin{equation*}
  \lim_{n\rightarrow+\infty} \phi(\xi_n)=(1+b)/2 \text{ and } \lim_{n\rightarrow+\infty} \phi'(\xi_n)=0.
\end{equation*}
Integrating $\phi$-equation of \eqref{eq:tr} from 0 to $\xi_n$, we obtain
\begin{equation}\label{eq:k1}
\begin{split}
  c[\phi(\xi_n)-\phi(0)]-&d_1\int_{0}^{\xi_n} \mathcal{N}_1\left[\phi\right](\xi)d\xi=\int_{0}^{\xi_n}\phi(\xi)f\left(\phi,\psi\right)(\xi)d\xi.
\end{split}
\end{equation}
Hence the left-hand of \eqref{eq:k1} is bounded similar to the proof in Lemma \ref{le:k1}, whereas the right-hand of \eqref{eq:k1} is unbounded since
\begin{equation}\label{eq:t3}
    \liminf\limits_{n\rightarrow+\infty}f\left(\phi,\psi\right)(\xi_n)\geq f\left(\phi_-,\psi_+\right)  \geq f\left((1+b)/2,1\right)>0
\end{equation}
under the assumption \eqref{eq:as1}, which leads to a contradiction. If $\phi(\xi)$ is oscillatory as $\xi\rightarrow +\infty$, we also can find a
sequence $\left\{\xi_n\right\}$ of the minimal points of $\phi(\xi)$ with $\xi_n\rightarrow+\infty$ as $n\rightarrow+\infty$ such that
\begin{equation*}
  \lim_{n\rightarrow+\infty}\phi(\xi_n)=(1+b)/2 \text{ and }\phi'(\xi_n)=0
\end{equation*}
Meanwhile, \eqref{eq:t3} and \eqref{eq:k1} still valid, which yields a contradiction. Similarly, we obtain $\phi_+<1$ by using blew inequality
\begin{align*}
\limsup\limits_{n\rightarrow+\infty}f\left(\phi,\psi\right)(\xi_n)\leq f\left(\phi_+,\psi_-\right)=f\left(1,\psi_-\right)=-\displaystyle\frac{m\psi_-}{1+a\psi_-}<0.
\end{align*}
Therefore, we complete the proof.

\qed

We are in this position to state the existence theorem for $c>c^*$.
\begin{theorem}\label{th:c1}
Under the assumption \eqref{eq:as1},
system \eqref{eq:rd} has traveling waves connecting $(1,0)$ and $(u_2^*,u_2^*)$ for $c>c^*$.
\end{theorem}
\Proof We define a sequence $\left\{\gamma_n\right\}_{n=-1}^{+\infty}$, where
\begin{align*}
\left\{
\begin{array}{ll}
\gamma_{-1}=\displaystyle\frac{1+b}{2},\ \ \gamma_{0}=1,\\[0.4cm]
\gamma_{n+1}=\displaystyle\frac{1}{2}\left(b+1+\sqrt{(1-b)^2-\frac{4mb\gamma_{n}}{\gamma_{n-1}+a\gamma_{n}}}\right).
\end{array}
\right.
\end{align*}
Under the assumption \eqref{eq:as1}, $\gamma_n$ is well-defined and $(1+b)/2\leq\gamma_n\leq1$ for all $n$. Moreover, note that
\begin{equation*}
(1-\gamma_{n+1})\left(\displaystyle\frac{\gamma_{n+1}}{b}-1\right)=\displaystyle\frac{m\gamma_{n}}{\gamma_{n-1}+a\gamma_{n}},
\end{equation*}
based on the fact that function $(1-\gamma)\left(\gamma/b-1\right)$ is decreasing in $(1+b)/2\leq\gamma\leq1$, we first show the following claims hold:\\

\textbf{Claim 1.} Under the assumption \eqref{eq:as1}, the sequences $\left\{\gamma_{2n}\right\}$ and $\left\{\gamma_{2n+1}\right\}$ are adjacent, that is
\begin{equation}\label{eq:px}
  \gamma_{-1}<\gamma_1<\cdot\cdot\cdot<\gamma_{2n-1}<\cdot\cdot\cdot<u_2^*<\cdot\cdot\cdot<\gamma_{2n}<\cdot\cdot\cdot<\gamma_{2}<\gamma_{0}.
\end{equation}
Note that $\gamma_{-1}<\gamma_{1}<u_2^*<\gamma_{2}<\gamma_{0}$ holds since
\begin{align*}
  &(1-\gamma_{1})\left(\displaystyle\frac{\gamma_{1}}{b}-1\right)=\displaystyle\frac{2m}{1+b+2a}>\frac{m}{1+a}=
  (1-u_2^*)\left(\displaystyle\frac{u_2^*}{b}-1\right),\\[0.2cm]
  &(1-\gamma_{1})\left(\displaystyle\frac{\gamma_{1}}{b}-1\right)=\displaystyle\frac{2m}{1+b+2a}<\frac{(1-b)^2}{4b}=
  (1-\gamma_{-1})\left(\displaystyle\frac{\gamma_{-1}}{b}-1\right),\\[0.2cm]
  &(1-\gamma_{2})\left(\displaystyle\frac{\gamma_{2}}{b}-1\right)=\displaystyle\frac{m\gamma_{1}}{\gamma_{0}+a\gamma_{1}}<\frac{m}{1+a}=
  (1-u_2^*)\left(\displaystyle\frac{u_2^*}{b}-1\right),\\[0.2cm]
  &(1-\gamma_{2})\left(\displaystyle\frac{\gamma_{2}}{b}-1\right)=\displaystyle\frac{m\gamma_{1}}{\gamma_{0}+a\gamma_{1}}>0=
  (1-\gamma_{0})\left(\displaystyle\frac{\gamma_{0}}{b}-1\right).
\end{align*}
We assume that
\begin{equation}\label{eq:z1}
  \gamma_{n-1}<\gamma_{n+1}<u_2^*<\gamma_{n+2}<\gamma_{n} \text{ with } n=2k, \ k\in\mathbb N,
\end{equation}
and will show that
\begin{equation}\label{eq:s}
  \gamma_{n+1}<\gamma_{n+3}<u_2^*<\gamma_{n+4}<\gamma_{n+2} \text{ with } n=2k,\ k\in\mathbb N.
\end{equation}
Due to $\gamma_{n+2}>u_2^*>\gamma_{n+1}$, we have
\begin{align*}
  &(1-\gamma_{n+3})\left(\displaystyle\frac{\gamma_{n+3}}{b}-1\right)=\displaystyle\frac{m\gamma_{n+2}}{\gamma_{n+1}+a\gamma_{n+2}}>\displaystyle\frac{m}{1+a}=(1-u_2^*)\left(\displaystyle\frac{u_2^*}{b}-1\right),
\end{align*}
hence $u_2^*>\gamma_{n+3}$.  Thanks to $\gamma_{n+2}>u_2^*>\gamma_{n+3}$, we further get
\begin{equation*}
  (1-\gamma_{n+4})\left(\displaystyle\frac{\gamma_{n+4}}{b}-1\right)=\displaystyle\frac{m\gamma_{n+3}}{\gamma_{n+2}+a\gamma_{n+3}}<\displaystyle\frac{m}{1+a}=(1-u_2^*)\left(\displaystyle\frac{u_2^*}{b}-1\right).
\end{equation*}
which implies that $\gamma_{n+4}>u_2^*$.
On the other hand, from \eqref{eq:z1}, we have
\begin{align*}
  (1-\gamma_{n+3})\left(\displaystyle\frac{\gamma_{n+3}}{b}-1\right)=\displaystyle\frac{m\gamma_{n+2}}{\gamma_{n+1}+a\gamma_{n+2}}<\displaystyle\frac{m\gamma_{n}}{\gamma_{n-1}+a\gamma_{n}}=(1-\gamma_{n+1})\left(\displaystyle\frac{\gamma_{n+1}}{b}-1\right),
\end{align*}
hence $\gamma_{n+3}>\gamma_{n+1}$. Due to $\gamma_{n+3}>\gamma_{n+1}$ and $\gamma_{n}>\gamma_{n+2}$, we further get
\begin{equation*}
  (1-\gamma_{n+4})\left(\displaystyle\frac{\gamma_{n+4}}{b}-1\right)=\displaystyle\frac{m\gamma_{n+3}}{\gamma_{n+2}+a\gamma_{n+3}}>\displaystyle\frac{m\gamma_{n+1}}{\gamma_{n}+a\gamma_{n+1}}=(1-\gamma_{n+2})\left(\displaystyle\frac{\gamma_{n+2}}{b}-1\right).
\end{equation*}
Hence we readily get \eqref{eq:s}, which further asserts \eqref{eq:px}. Therefore, the claim is valid by mathematical induction method.\\

\textbf{Claim 2.} The sequences $\gamma_{2n}$ and $\gamma_{2n-1}$ converge to $u_2^*$, respectively.\\

From \textbf{Claim 1}, the sequence $\gamma_{2n}$ is decreasing with $u_2^*$ as the lower bound and the sequence $\gamma_{2n-1}$ is increasing with $u_2^*$ as the upper bound, then there are two constants $\gamma^*$ and  $\gamma_*$ with $\gamma^*\geq u_2^*\geq \gamma_*$ such that $\gamma_{2n}\rightarrow \gamma^*$ and $\gamma_{2n-1}\rightarrow \gamma_*$ as $n\rightarrow+\infty$. Note that the inequality below
\begin{equation*}
  |\gamma^*-\gamma_*|\leq|\gamma^*-\gamma_{2n}|+|\gamma_{2n}-\gamma_{2n-1}|+|\gamma_{2n-1}-\gamma_*|,
\end{equation*}
therefore, it is sufficient to prove that the sequences $\gamma_{2n}-\gamma_{2n-1}\rightarrow0$ as $n\rightarrow+\infty$. Note that
\begin{equation*}
  \max\limits_{(1+b)/2\leq x,y\leq1}\left\{\frac{\partial}{\partial x}\left(\frac{x}{y+ax}\right), \frac{\partial}{\partial y}\left(\frac{x}{y+ax}\right)\right\}\leq\frac{4}{(1+a)^2(1+b)^2},
\end{equation*}
by using mean value theorem, we have for all $n\in\mathbb N^+$
\begin{align*}
  &\frac{b\gamma_{2n-2}}{\gamma_{2n-3}+a\gamma_{2n-2}}-\frac{b\gamma_{2n-1}}{\gamma_{2n-2}+a\gamma_{2n-1}}\\[0.2cm]
\leq&\frac{4b}{(1+a)^2(1+b)^2}\left[(\gamma_{2n-2}-\gamma_{2n-3})+(\gamma_{2n-2}-\gamma_{2n-1})\right],\\[0.2cm]
<&\frac{8b}{(1+a)^2(1+b)^2}(\gamma_{2n-2}-\gamma_{2n-3}),\\[0.2cm]
<&\frac{2}{(1+a)^2}(\gamma_{2n-2}-\gamma_{2n-3})
\end{align*}
owing to $\gamma_{2n-3}<\gamma_{2n-1}$. Then a direct calculation gives
\begin{align*}
  \gamma_{2n}-\gamma_{2n-1}&=\frac{1}{2}\left(\sqrt{(1-b)^2-\frac{4mb\gamma_{2n-1}}{\gamma_{2n-2}+a\gamma_{2n-1}}}-\sqrt{(1-b)^2-\frac{4mb\gamma_{2n-2}}{\gamma_{2n-3}+a\gamma_{2n-2}}}\right)\\[0.2cm]
  &=\frac{2m\left(\frac{b\gamma_{2n-2}}{\gamma_{2n-3}+a\gamma_{2n-2}}-\frac{b\gamma_{2n-1}}{\gamma_{2n-2}+a\gamma_{2n-1}}\right)}{\sqrt{(1-b)^2-\frac{4mb\gamma_{2n-1}}{\gamma_{2n-2}+a\gamma_{n2-1}}}+\sqrt{(1-b)^2-\frac{4mb\gamma_{2n-2}}{\gamma_{2n-3}+a\gamma_{2n-2}}}}\\[0.2cm]
  &\leq\frac{2m\left(\frac{b\gamma_{2n-2}}{\gamma_{2n-3}+a\gamma_{2n-2}}-\frac{b\gamma_{2n-1}}{\gamma_{2n-2}+a\gamma_{2n-1}}\right)}{\sqrt{(1-b)^2-\frac{4mb\gamma_{2n-1}}{\gamma_{2n-2}+a\gamma_{2n-1}}}}\\[0.2cm]
  &\leq\frac{ \frac{4m}{(1+a)^2}}{\sqrt{(1-b)^2-\frac{4mb}{1+a}}}(\gamma_{2n-2}-\gamma_{2n-3})\\[0.2cm]
  &:=\rho(\gamma_{2n-2}-\gamma_{2n-3}).
\end{align*}
If $\rho\in(0,1)$, that is
\begin{equation*}
 m<\frac{(1+a)^3}{8}\left(\sqrt{b^2+4\left(\frac{1-b}{1+a}\right)^2}-b\right),
\end{equation*}
one can lightly verify the sequences $\gamma_{2n}-\gamma_{2n-1}\rightarrow0$ as $n\rightarrow+\infty$. Hence, the claim is valid.\\

Now, we show that $(\phi,\psi)(+\infty)=(u_2^*,u_2^*)$. From Lemma \ref{le:0}, we have
\begin{equation*}
\gamma_{-1}<\phi_-\leq\psi_-\leq\psi_+\leq\phi_+<\gamma_{0}.
\end{equation*}
From \textbf{Claim 2}, knowing that $(\phi,\psi)(+\infty)=(u_2^*,u_2^*)$ holds as long as for all $n\in\mathbb N^+$
\begin{equation*}
  \gamma_{2n-1}<\phi_-\leq\psi_-\leq\psi_+\leq\phi_+<\gamma_{2n}.
\end{equation*}
Hence we define
\begin{equation*}
  n_0=\sup\left\{n\in\mathbb N^+|\ \gamma_{2n-1}<\phi_-\leq\psi_-\leq\psi_+\leq\phi_+<\gamma_{2n}\right\},
\end{equation*}
and we will show that $n_0=+\infty$. By contradiction, we assume that $n_0$ is finite.
From the definition of $n_0$, it holds that either
\begin{equation*}
  \gamma_{2n_0-1}=\phi_-\leq\psi_-\leq\psi_+\leq\phi_+\leq\gamma_{2n}
\end{equation*}
or
\begin{equation*}
  \gamma_{2n_0-1}\leq\phi_-\leq\psi_-\leq\psi_+\leq\phi_+=\gamma_{2n}.
\end{equation*}
Firstly, we assume that $\phi_-=\gamma_{2n_0-1}$. If $\phi(\xi)$ is eventually monotone, then $\phi(+\infty)=\gamma_{2n_0-1}$ and
\begin{equation*}
  \int_{0}^{+\infty} \phi'(\xi)d\xi=\gamma_{2n_0-1}-\phi(0)<\infty
\end{equation*}
since $\phi(\xi)$ is bounded on $\mathbb{R}$.
Furthermore, we have
$\limsup_{\xi\rightarrow+\infty}\phi'(\xi)=0$ if  $\phi'(\xi)\leq0$  for $\xi\gg 1$ or
$\liminf_{\xi\rightarrow+\infty}\phi'(\xi)=0$ if  $\phi'(\xi)\geq0$  for $\xi\gg 1$.
Then we can find a sequence $\left\{\xi_n\right\}$ with ${\xi_n\rightarrow+\infty}$ as $n\rightarrow+\infty$ such that
\begin{equation*}
  \lim_{n\rightarrow+\infty} \phi(\xi_n)=\gamma_{2n_0-1} \text{ and } \lim_{n\rightarrow+\infty} \phi'(\xi_n)=0.
\end{equation*}
Integrating $\phi$-equation of \eqref{eq:tr} from 0 to $\xi_n$,  we obtain
\begin{equation}\label{eq:k3}
\begin{split}
  c[\phi(\xi_n)-\phi(0)]-&d_1\int_{0}^{\xi_n} \mathcal{N}_1\left[\phi\right](\xi)d\xi=\int_{0}^{\xi_n}\phi(\xi)f\left(\phi,\psi\right)(\xi)d\xi.
\end{split}
\end{equation}
Hence the left-hand of \eqref{eq:k3} is bounded similar to the proof in Lemma \ref{le:k1}, whereas the right-hand of \eqref{eq:k3} is unbounded since
\begin{equation}\label{eq:t4}
\begin{split}
  \liminf\limits_{n\rightarrow+\infty}f\left(\phi,\psi\right)(\xi_n)&\geq f\left(\gamma_{2n_0-1},\gamma_{2n_0-2}\right)\\[0.2cm]
  &>(1-\gamma_{2n_0-1})\left(\frac{\gamma_{2n_0-1}}{b}-1\right)-\displaystyle\frac{m\gamma_{2n_0-2}}{\gamma_{2n_0-3}+a\gamma_{2n_0-2}}=0,
\end{split}
\end{equation}
which leads a contradiction. If $\phi(\xi)$ is oscillatory as $\xi\rightarrow +\infty$, we can find a sequence $\left\{\xi_n\right\}$ of the minimal points of $\phi(\xi)$ with $\xi_n\rightarrow+\infty$ as $n\rightarrow+\infty$ such that
\begin{equation*}
  \lim_{n\rightarrow+\infty}\phi(\xi_n)=\gamma_{2n_0-1} \text{ and }\phi'(\xi_n)=0.
\end{equation*}
We further get \eqref{eq:k3}, which contradicts to \eqref{eq:t4} again. The case $\phi_+=\gamma_{2n_0}$ can be treated similarly, since
\begin{equation*}
\begin{split}
  \limsup\limits_{n\rightarrow+\infty}f\left(\phi,\psi\right)(\xi_n)&\leq f\left(\gamma_{2n_0},\gamma_{2n_0-1}\right)\\[0.2cm]
  &<(1-\gamma_{2n_0})\left(\frac{\gamma_{2n_0}}{b}-1\right)-\displaystyle\frac{m\gamma_{2n_0-1}}{\gamma_{2n_0-2}+a\gamma_{2n_0-1}}=0.
\end{split}
\end{equation*}
Hence we must have $n_0=+\infty$, which implies that $(\phi,\psi)(+\infty)=(u_2^*,u_2^*)$ and finishes this proof.

\qed

\section{The existence of traveling waves for $c=c^*$}
\noindent

This section is devoted to the existence of traveling waves to system \eqref{eq:rd} for $c=c^*$, which depends on a limiting argument. Let us start with the following result.

\begin{lemma}
For $c>c^*$, the solution $(\phi,\psi)(\xi)$ of system \eqref{eq:tr}-\eqref{eq:b} satisfies
\begin{equation*}
  \lim\limits_{\xi\rightarrow-\infty}\frac{\psi'(\xi)}{\psi(\xi)}=\lambda\in\left\{\lambda_1,\ \lambda_2\right\}.
\end{equation*}
\end{lemma}
\Proof From Theorem \ref{th:str}, we have $\psi(\xi)>0$ over $\mathbb R$. Let
\begin{equation*}
  Z(\xi)=\frac{\psi'(\xi)}{\psi(\xi)} \text{ and } B(\xi)=g(\xi)-d_2.
\end{equation*}
then  $\psi(\xi)$-equation gives
\begin{equation*}
  cZ(\xi)=d_2\int_{\mathbb R} J_2(y)e^{\int_{\xi}^{\xi-y}Z(s)ds}dy+B(\xi).
\end{equation*}
Since $(\phi,\psi)(-\infty)=(1,0)$, we have
\begin{equation*}
  B(-\infty)=g(-\infty)-d_2=s\left(1-\frac{\psi(-\infty)}{\phi(-\infty)}\right)-d_2=s-d_2.
\end{equation*}
By Lemma \ref{pr:1}, we obtain $\lambda:= \lim\limits_{\xi\rightarrow-\infty}Z(\xi)$ exist, and satisfies
$$c\lambda=d_2\left(\int_{\mathbb R} J_2(y)e^{-\lambda y}dy-1\right)+s.$$
According to Lemma \ref{le:cx}, we complete the proof.

\qed

Now, we choose a strictly decreasing sequence $\left\{c_n\right\}$ with $c_n\in(c^*, c^*+ 1)$ and $\lim_{n\rightarrow+\infty}c_n=c^*$. Thus, for each $c_n$, there exists a positive solution $(\phi_n,\psi_n)(\xi)$ of system \eqref{eq:tr}-\eqref{eq:b}.
From the above lemma, we have $\psi_n'(\xi)>0$ for $\xi\ll-1$ since $\psi_n(\xi)>0$ over $\mathbb R$. On the other hand, since $(\phi_n,\psi_n)(\xi+ a)$ for any $a\in\mathbb R$ is also solution of system \eqref{eq:tr}-\eqref{eq:b}. Therefore, we assume that there exists a positive constant $\delta<\min\left\{(1+b)/8, u_1^*/2\right\}$ such that $\psi_n(0)=\delta$ and $\psi_n(\xi)\leq\delta$ for $\xi<0$.

\begin{theorem}
Under the assumption \eqref{eq:as1}, system \eqref{eq:rd} has a traveling wave connecting $(1,0)$ and $(u_2^*,u_2^*)$ for $c=c^*$.
\end{theorem}
\textbf{Proof.} From Theorem \ref{th:str}, we have $(1+b)/2<\phi_{n}(\xi)<1$ and $0<\psi_{n}(\xi)<1$ over $\mathbb R$, then $\phi'_{n}(\xi)$ and $\psi'_{n}(\xi)$ is uniformly bounded over $\mathbb R$ from system \eqref{eq:tr}, which implies that $\phi_{n}(\xi)$ and $\psi_{n}(\xi)$ are equicontinuous.  Additionally, it follows from \emph{(J2)} that
\begin{align*}
  \left|\frac{d}{d \xi} \int_{\mathbb{R}} J_{1}(\xi-y) \phi_n(y) \mathrm{d} y\right|=\left|\int_{\mathbb{R}} \frac{d}{d \xi} J_{1}(\xi-y) \phi_n(y) \mathrm{d} y\right| \leqslant \int_{\mathbb{R}}\left|J_{1}^{\prime}(y)\right| d y,\\[0.2cm]
  \left|\frac{d}{d \xi} \int_{\mathbb{R}} J_{2}(\xi-y) \psi_n(y) \mathrm{d} y\right|=\left|\int_{\mathbb{R}} \frac{d}{d \xi} J_{2}(\xi-y) \psi_n(y) \mathrm{d} y\right| \leqslant \int_{\mathbb{R}}\left|J_{2}^{\prime}(y)\right| d y .
\end{align*}
Then by calculating derivative on $\xi$ in system \eqref{eq:tr}, we know that $\phi''_{n}(\xi)$ and $\psi''_{n}(\xi)$ are uniformly bounded over $\mathbb R$, which implies that $\phi'_{n}(\xi)$ and $\psi'_{n}(\xi)$ are equicontinuous. Hence, by Arzela-Ascoli theorem, up to extracting a subsequence for necessary, there are functions $\phi(\xi)$ and $\psi(\xi)$, such that $(\phi_{n},\psi_{n})(\xi)$ and $(\phi'_{n},\psi'_{n})(\xi)$ converge uniformly to $(\phi, \psi)(\xi)$ and $(\phi', \psi')(\xi)$
on every bounded interval and point-wise over $\mathbb R$. Furthermore, the dominated convergence theorem yields that
\begin{align*}
&\int_{\mathbb{R}} J_{1}(\xi-y) \phi_n(y) \mathrm{d} y \rightarrow \int_{\mathbb{R}} J_{1}(\xi-y) \phi(y) \mathrm{d} y  \text{ as } n \rightarrow+\infty, \\[0.2cm]
&\int_{\mathbb{R}} J_{2}(\xi-y) \psi_n(y) \mathrm{d} y \rightarrow \int_{\mathbb{R}} J_{2}(\xi-y) \psi(y) \mathrm{d} y  \text{ as } n \rightarrow+\infty.
\end{align*}
Therefore, we know immediately that $(\phi, \psi)(\xi)$ is a solution of system \eqref{eq:tr} with $c=c^*$ by letting $n\rightarrow+\infty$ in system \eqref{eq:tr} with $(\phi, \psi)(\xi)=(\phi_n, \psi_n)(\xi)$. Simultaneously, we also have $(1+b)/2\leq\phi(\xi)\leq1$ and $0\leq\psi(\xi)\leq1$ over $\mathbb R$. Similar to proof in Theorem \ref{th:str}, we obtain that $(1+b)/2<\phi(\xi)<1$ and $0<\psi(\xi)<1$ over $\mathbb R$.

It should be emphasized that the proofs of Lemmas \ref{le:k1} and \ref{le:0}, and Theorem \ref{th:c1} are independence of $c$, then we still get $(\phi, \psi)(+\infty)=(u_2^*,u_2^*)$.
Therefore, what we need to do is to prove $(\phi, \psi)(-\infty)=(1,0)$. Naturally, we define
\begin{align*}
&\overline{\phi}=\limsup \limits_{\xi \rightarrow-\infty} \phi(\xi),\quad \underline{\phi}=\liminf\limits_{\xi \rightarrow-\infty} \phi(\xi), \\[0.2cm]
&\overline{\psi}=\limsup \limits_{\xi \rightarrow-\infty} \psi(\xi),\quad \underline{\psi}=\liminf\limits_{\xi \rightarrow-\infty} \psi(\xi),
\end{align*}
and then divide our proof into the following two cases.

\textbf{Case 1.} $\overline{\phi}=\underline{\phi}$. In this case, we will prove that $\psi(-\infty)$ exists and equals to $0$. Actually, assume that $\underline{\psi}<\overline{\psi}$, there exist two sequences $\left\{x_n\right\}$ and $\left\{y_n\right\}$ satisfying $x_n$, $y_n\rightarrow-\infty$ as $n\rightarrow+\infty$ such that
\begin{equation*}
  \lim\limits_{n\rightarrow+\infty}\psi(x_n)=\underline{\psi} \text{ and }   \lim\limits_{n\rightarrow+\infty}\psi(y_n)=\overline{\psi}.
\end{equation*}
From Lemma \ref{pr:2}, we have $\phi'(-\infty)=0$ since $\phi(-\infty)$ exists.
Using Lemma \ref{pr:3}, for any sequence $\left\{\tau_n\right\}$ with $\tau_n\rightarrow-\infty$ as $n\rightarrow+\infty$, we get
\begin{equation*}
  \lim\limits_{n\rightarrow+\infty}\mathcal{N}_1[\phi](\tau_n)
  =0.
\end{equation*}
Selecting $\left\{\tau_n\right\}$ as $\left\{x_n\right\}$ and $\left\{y_n\right\}$ respectively, and letting $n\rightarrow+\infty$, we arrive at
\begin{align*}
  \left\{
  \begin{array}{c}
  \phi(-\infty)f\left(\phi(-\infty),\underline{\psi}\right)=0,\\[0.2cm]
  \phi(-\infty)f\left(\phi(-\infty),\overline{\psi}\right)=0.
  \end{array}
 \right.
\end{align*}
Since $\phi(\xi)\geq(1+b)/2$ over $\mathbb R$, we must have
\begin{equation*}
  f\left(\phi(-\infty),\underline{\psi}\right)=f\left(\phi(-\infty),\overline{\psi}\right)=0.
\end{equation*}
From the expression of $f(\phi,\psi)$, automatically,
\begin{equation*}
  \displaystyle\frac{m\underline{\psi}}{\phi(-\infty)+a\underline{\psi}}=\displaystyle\frac{m\overline{\psi}}{\phi(-\infty)+a\overline{\psi}},
\end{equation*}
which yields that $\underline{\psi}=\overline{\psi}$ and $\psi(-\infty)$ exists. From Lemma \ref{pr:2}, we have $\psi'(-\infty)=0$. Similarly, for any sequence $\left\{\tau_n\right\}$ with $\tau_n\rightarrow-\infty$ as $n\rightarrow+\infty$,
\begin{equation*}
  \lim\limits_{n\rightarrow+\infty}\mathcal{N}_2[\psi](\tau_n)=0.
\end{equation*}
Hence, we have
\begin{equation*}
  \left\{
  \begin{split}
  &\phi(-\infty)f\left(\phi(-\infty),\psi(-\infty)\right)=0,\\[0.2cm]
  &\psi(-\infty)g\left(\phi(-\infty),\psi(-\infty)\right)=0.
  \end{split}
 \right.
\end{equation*}
Thus, one of four cases may happen:
\begin{align*}
  &\left(\phi(-\infty),\psi(-\infty)\right)=(b,0),\\[0.2cm]
  &\left(\phi(-\infty),\psi(-\infty)\right)=(1,0),\\[0.2cm]
  &\left(\phi(-\infty),\psi(-\infty)\right)=(u_1^*,u_1^*),\\[0.2cm]
  &\left(\phi(-\infty),\psi(-\infty)\right)=(u_2^*,u_2^*).
\end{align*}
Since $\delta<u_1^*/2$ and $\phi(\xi)\geq(1+b)/2$, we have $(\phi,\psi)(-\infty)=(1,0)$.\\

\textbf{Case 2.} $\overline{\phi}\neq\underline{\phi}$. We claim that $\psi(-\infty)=0$. Otherwise, there exist $\zeta\in(0,\delta)$ such that
\begin{equation}\label{eq:d}
  \limsup\limits_{\xi\rightarrow-\infty}\psi(\xi)=\zeta.
\end{equation}
Up to extracting a subsequence, there is a sequence $\left\{\ell_n\right\}$ with $\ell_n\rightarrow-\infty$ as $n\rightarrow+\infty$ such that
\begin{equation*}
    \psi(\ell_n)>\frac{\zeta}{2} \text{ for all } n.
\end{equation*}
From the uniform continuity of $\psi(\xi)$, we have for an appropriate $\epsilon_1>0$
\begin{equation*}
   \psi(\ell_n+\xi)>\frac{\zeta}{4} \text{ for } \xi\in(-\epsilon_1,\epsilon_1).
\end{equation*}
Now, we consider initial value problem
\begin{equation*}
  \left\{
  \begin{split}
    &\varphi_t= d_2\mathcal{N}_2[\varphi](x,t)+s\varphi\left(1-\frac{2\varphi}{1+b}\right), \\[0.2cm]
    &\varphi(x,0)=\varphi(x),
  \end{split}
  \right.
\end{equation*}
where $\varphi(x)$ satisfies the following conditions:
\begin{description}
  \item (1) $\varphi(x)$ is uniformly continuous on $x$,
  \item (2) $\varphi(x)=\zeta/4$ for $x\in[-\epsilon_1/2,\epsilon_1/2]$,
  \item (3) $\varphi(x)$ is decreasing for $x\in[\epsilon_1/2,\epsilon_1]$ and increasing for $x\in[-\epsilon_1/2,-\epsilon_1]$,
  \item (4) $\varphi(x)=0$ for $|x|>\epsilon_1$.
\end{description}
Spreading speed theory \cite{zy22} gives for  $c\in(0,c^*)$
\begin{equation}\label{eq:f1}
  \liminf\limits_{t\rightarrow+\infty}\inf\limits_{|x|<ct}\varphi(x,t)=\frac{1+b}{2}.
\end{equation}
Thanks to $\delta<(1+b)/8$, there is a constant $T>0$ such that for $c\in(0,c^*)$
\begin{equation*}
  \inf\limits_{|x|<ct}\varphi(x,t)>2\delta.
\end{equation*}
For the above $T$, we further choose two subsequences $\left\{\ell_{1n}\right\}$ and $\left\{\ell_{2n}\right\}$ satisfying  for all $n$
\begin{equation*}
  \ell_{1n}-\ell_{2n}>c^*T,\ \psi(\ell_{1n})>\frac{\zeta}{2} \text{ and } \psi(\ell_{2n})>\frac{\zeta}{2}.
\end{equation*}
From $\psi$-equation, the function $w(x,t):=\psi(x+c^*t+\ell_{2n})$ satisfies
\begin{equation*}
  \left\{
  \begin{split}
    &w_t\geq d_2\mathcal{N}_2[w](x,t)+sw\left(1-\frac{2w}{1+b}\right), \\[0.2cm]
    &w(x,0)=\psi(x+\ell_{2n}).
  \end{split}
  \right.
\end{equation*}
By using comparison principle \cite{zy22}, we have for $c\in(0,c^*)$
\begin{equation}\label{eq:ad}
  \liminf\limits_{t\rightarrow+\infty}\inf\limits_{|x|<ct}w(x,t)\geq\liminf\limits_{t\rightarrow+\infty}\inf\limits_{|x|<ct}\varphi(x,t)=(1+b)/2.
\end{equation}
Now, we fix $x=0$ and $t=(\ell_{1n}-\ell_{2n})/c^*$. Obviously, $|x|<ct$ for $c\in(0,c^*)$. Then for $t>T$
\begin{equation*}
  w(0,t)=w(0,(\ell_{1n}-\ell_{2n})/c^*)=\psi(\ell_{1n})>\varphi(0,t)>2\delta>\zeta
\end{equation*}
by the choice of $\zeta$. Hence we obtain
\begin{equation*}
  \limsup\limits_{\xi\rightarrow-\infty}\psi(\xi)>\zeta,
\end{equation*}
which contradicts to \eqref{eq:d}. Therefore, we have $\psi(-\infty)=0$. Under the condition $\overline{\phi}\neq\underline{\phi}$, there are sequences $\left\{x_n\right\}$ and $\left\{y_n\right\}$ with $x_n$, $y_n\rightarrow-\infty$ as $n\rightarrow+\infty$ such that
\begin{align*}
  \lim\limits_{n\rightarrow+\infty}\phi(x_n)=\underline{\phi} \text{ and } \phi'(x_n)=0,\\[0.2cm]   \lim\limits_{n\rightarrow+\infty}\phi(y_n)=\overline{\phi}  \text{ and } \phi'(y_n)=0.
\end{align*}
From Lemma \ref{pr:3}, taking $\xi=x_n$ or $\xi=y_n$ in $\phi$-equation, and letting $n\rightarrow+\infty$, we have
\begin{align*}
  \underline{\phi}\left(1-\underline{\phi}\right)\left(\frac{\underline{\phi}}{b}-1\right)\leq0,\\[0.2cm]
  \overline{\phi}\left(1-\overline{\phi}\right)\left(\frac{\overline{\phi}}{b}-1\right)\geq0.
\end{align*}
Since $\overline{\phi}>\underline{\phi}\geq(1+b)/2$, we have $\underline{\phi}<\overline{\phi}\leq1\leq\underline{\phi}$, which suggests that $\overline{\phi}\neq\underline{\phi}$ cannot occur.\\

Therefore, we complete the proof by combining \textbf{Case 1} with \textbf{Case 2}.

\qed

Finally, let us finish the proof of Theorem \ref{th:ma} by proving that traveling waves does not exists for $0<c<c^*$.

\begin{theorem}
For $0<c<c^*$, system \eqref{eq:rd} has no traveling waves connecting $(1,0)$ and $(u_2^*,u_2^*)$.
\end{theorem}
\textbf{Proof.} By contradiction, we assume that there exists a positive solution $(\phi, \psi)(\xi)$ of \eqref{eq:tr}-\eqref{eq:b} for given $0<c_0<c^*$. Through boundary conditions \eqref{eq:b}, we have $\phi(\xi)\geq1/K$ over $\mathbb R$ for some $K\gg1$.
Let $w(x,t):=\psi(x+c_0t)$, from $\psi$-equation, we get
\begin{equation*}
\left\{\begin{split}
&w_t\geq d_{2} \mathcal{N}_2[w](x,t)+s w(x,t)\left(1-Kw(x,t)\right),\\[0.2cm]
&w(0,x)=\psi(x).
\end{split}
\right.
\end{equation*}
Spreading speed theory and comparison principle \cite{zy22} give
\begin{equation*}
  \lim\limits_{t\rightarrow+\infty}\inf\limits_{2|x|=(c_0+c^*)t} w(x,t)\geq \frac{1}{K}.
\end{equation*}
As $2x=-(c_0+c^*)t$, we have
\begin{equation*}
  w(x,t)=\psi(x+c_0t)=\psi\left(-\frac{(c_0+c^*)t}{2}+c_0t\right)=\psi\left(\frac{(c_0-c^*)t}{2}\right)\geq \frac{1}{K},
\end{equation*}
which contradicts to $\psi(-\infty)=0$. Therefore, we complete the proof.

\qed

\setcounter{lemma}{0}
\renewcommand{\thelemma}{A.\arabic{lemma}}
\section*{Appendix}
In appendix, some important lemmas used in this paper are given.
\begin{lemma}\label{pr:1}
(\cite{zlw12} ) Assume $c>0$ and $B(\xi)$ is a continuous function with $B(\infty) := \lim_{\xi\rightarrow\infty}B(\xi)$. Let $Z(\xi)$ be a measurable function satisfying
\begin{equation*}
  cZ(\xi)=\int_{\mathbb R} J_i(y)e^{\int_{\xi}^{\xi-y}Z(s)ds}dy+B(\xi) \text{ in } \mathbb R,\ \ i=1,2.
\end{equation*}
Then, $Z(\xi)$ is uniformly continuous and bounded. Moreover, $\mu:= \lim_{\xi\rightarrow\infty}Z(\xi)$ exist and are real roots of the characteristic equation
\begin{equation*}
  c\mu=\int_{\mathbb R} J_i(y)e^{-\mu y}dy+B(\infty),\ \ i=1,2.
\end{equation*}
\end{lemma}
\begin{lemma}\label{pr:2}
(\cite{b59} ) Assume that $w(\xi)\in C^1(b, +\infty)$ and $\lim_{\xi\rightarrow+\infty}w(\xi)$ exists. If $w'(\xi)$ is uniformly continuous, then $\lim_{\xi\rightarrow+\infty}w'(\xi)=0$.
\end{lemma}
\begin{lemma}\label{pr:3}
(\cite{w21} ) Assume that $J(\xi)\geq0$ and $\displaystyle\int_{\mathbb R}J(\xi)d\xi=1$, and $\omega(\xi)$ is a nonnegative bounded continuous function on $\mathbb R$. Then we have
\begin{equation*}
\begin{split}
&\liminf _{\xi \rightarrow \infty} \int_{\mathbb{R}} J(y) \omega(\xi-y) d y \geq \liminf _{\xi \rightarrow \infty} \omega(\xi):=\omega^-, \\[0.2cm]
&\limsup _{\xi \rightarrow \infty} \int_{\mathbb{R}} J(y) \omega(\xi-y) d y \leq \limsup _{\xi \rightarrow \infty} \omega(\xi):=\omega^+.
\end{split}
\end{equation*}
In particular, if $\omega(\infty)$ exists, that is, $\omega^-=\omega^+=\omega(\infty)$, then
\begin{equation*}
  \lim _{\xi \rightarrow \infty} \int_{\mathbb{R}} J(y) \omega(\xi-y) d y=\omega(\infty).
\end{equation*}
\end{lemma}
\section*{Acknowledgments}
This work is supported by the National Natural Science Foundation of China (No. 12171039 and 12271044).


\bibliographystyle{elsarticle-num}

\end{document}